\documentclass[a4paper,12pt]{article}
\catcode`\"=13          
\def"#1{\v#1}
\usepackage[dvips]{epsfig,graphicx}
\usepackage{amsmath,amsfonts,amsbsy} 

\newcommand{\CC}{{\mathbb C}}

\newcommand{\G}{\Gamma}
\newcommand{\ov}{\overline}
\newcommand{\la}{\langle}
\newcommand{\h}{\hat}
\newcommand{\ra}{\rangle}

\newcommand{\As}{A^*}
\newcommand{\Es}{E^*}
\newcommand{\mat}{{\rm Mat}}
\newcommand{\qed}{\hfill\hbox{\rule{3pt}{6pt}}}
\newcommand{\proof}{{\sc Proof. }}
\newcommand{\MX}{\mat_X(\CC)}

\newtheorem{theorem}{Theorem}[section]
\newtheorem{lemma}[theorem]{Lemma}
\newtheorem{corollary}[theorem]{Corollary}
\newtheorem{proposition}[theorem]{Proposition}
\newtheorem{definition}[theorem]{Definition}
\newtheorem{notation}[theorem]{Notation}

\title{THE TERWILLIGER ALGEBRA OF A DISTANCE-REGULAR GRAPH OF NEGATIVE TYPE}

\author{"Stefko Miklavi"c\footnote{Supported in part by ``Javna agencija za raziskovalno 
                                   dejavnost Republike Slovenije'', program no. Z1-9614.}
 \\ Primorska Institute for Natural Science and Technology  \\
        University of Primorska \\ 6000 Koper, Slovenia \\
        stefko.miklavic@upr.si}

\begin{document}
\maketitle
\begin{abstract}
  Let $\Gamma$ denote a distance-regular graph with diameter $D \ge 3$. 
  Assume $\Gamma$ has classical parameters $(D,b,\alpha,\beta)$ with $b < -1$. 
  Let $X$ denote the vertex set of $\G$ and let $A \in \MX$ denote the 
  adjacency matrix of $\G$. Fix $x \in X$ and let $\As \in \MX$ denote the corresponding dual 
  adjacency matrix. Let $T$ denote the subalgebra of $\MX$ generated by $A, \As$. We call $T$ the 
  {\em Terwilliger algebra} of $\G$ with respect to $x$. We show that 
  up to isomorphism there exist exactly two irreducible $T$-modules with endpoint $1$;
  their dimensions are $D$ and $2D-2$. For these $T$-modules we display a basis consisting
  of eigenvectors for $\As$, and for each basis we give the action of $A$.
\end{abstract}

\section{Introduction}
\label{sec:intro}

Let $\G$ denote a $Q$-polynomial distance-regular graph with diameter $D \ge 3$ and 
intersection numbers $a_i, b_i, c_i$ (see Section \ref{sec:prelim} for formal definitions). We recall 
the Terwilliger algebra of $\G$. Let $X$ denote the vertex set of $\G$ and let $A \in \MX$ 
denote the adjacency matrix of $\G$. Fix a ``base vertex'' $x \in X$ and let $\As \in \MX$ denote 
the corresponding dual adjacency matrix. Let $T=T(x)$ denote the subalgebra of $\MX$ generated by 
$A, \As$. The algebra $T$ is called the {\em Terwilliger algebra} of $\G$ with respect 
to $x$ \cite{ter1}. $T$ is closed under the conjugate-transpose map so $T$ is semi-simple 
\cite[Lemma 3.4(i)]{ter1}. Therefore each $T$-module is a direct sum of irreducible $T$-modules. 
Describing the irreducible $T$-modules is an active area of research [3--17], \cite{Mi1, Ta, ter1, TY}.

In this description there is an important parameter called the {\em endpoint} which we now recall.
Let $W$ denote an irreducible $T$-module. By the {\it endpoint} of $W$ we mean 
$\min \{i \mid 0 \le i \le D, \; \Es_i W \ne 0\}$, where $\Es_i \in \MX$ is the projection onto the 
$i$th subconstituent of $\G$ with respect to $x$ \cite[p. 378]{ter1}. There exists a unique 
irreducible $T$-module with endpoint $0$ \cite[Proposition 8.4]{egge1}; for a detailed description see 
\cite{Cu1, egge1}. 

Consider now the irreducible $T$-modules with endpoint $1$.
If $\G$ is bipartite, then these $T$-modules are described in 
\cite{Cu1,Cu2}. If $\G$ is nonbipartite with $a_1=0$, then these $T$-modules
are described in \cite{CMT,Mi1}. For the rest of this Introduction assume $a_1 \ne 0$. 
Assume further that $\G$ is of negative type and not a near polygon.
In \cite{Mi2} we described the combinatorial structure of $\G$. 
In the present paper we use this description to obtain the irreducible $T$-modules
that have endpoint 1. To summarize our results we note the following.
Let $W$ denote an irreducible 
$T$-module with endpoint $1$. Observe that $\Es_1 W$ is a $1$-dimensional eigenspace 
for $\Es_1 A \Es_1$ \cite[Theorem 2.2]{HI}. The corresponding
eigenvalue is called the {\em local eigenvalue} of $W$.
We show that up to isomorphism there exist exactly two irreducible $T$-modules with 
endpoint $1$. The first one has dimension $D$ and local eigenvalue $-1$. 
The second one has dimension $2D-2$ and local eigenvalue $a_1$. 
For these modules we display a basis consisting of eigenvectors for $\As$, 
and for each basis we give the action of $A$.
At present there is no classification of graphs that satisfy our assumptions;
see \cite[Section 6]{Mi2} for a summary of what is known. 
\section{Preliminaries}
\label{sec:prelim}

In this section we review some definitions and basic results concerning distance-regular graphs. 
See the book of Brouwer, Cohen and Neumaier \cite{BCN} for more background information.

\medskip \noindent
Let $\CC$ denote the complex number field and let $X$ denote a nonempty finite set. Let $\MX$ denote 
the $\CC$-algebra consisting of all matrices whose rows and columns are indexed by $X$ and whose 
entries are in $\CC$. Let $V=\CC^X$ denote the vector space over $\CC$ consisting of column vectors 
whose coordinates are indexed by $X$ and whose entries are in $\CC$. We observe $\MX$ acts on $V$ by 
left multiplication. We call $V$ the {\it standard module}. We endow $V$ with the Hermitean inner 
product $\langle \, , \, \rangle$ that satisfies $\langle u,v \rangle = u^t\overline{v}$ for 
$u,v \in V$, where $t$ denotes transpose and $\overline{\phantom{v}}$ denotes complex conjugation. For 
$y \in X$ let $\hat{y}$ denote the element of $V$ with a 1 in the $y$ coordinate and 0 in all 
other coordinates. We observe $\{\hat{y}\;|\;y \in X\}$ is an orthonormal basis for $V$.
The following will be useful: for each $B \in \MX$ we have 
\begin{equation}
\label{cez}
\la u, B v \ra = \la \ov{B}^t u, v \ra \qquad \qquad (u,v \in V).
\end{equation}

\medskip \noindent
Let $\G = (X,R)$ denote a finite, undirected, connected graph, without loops or multiple edges, 
with vertex set $X$ and edge set $R$. Let $\partial$ denote the path-length distance function for 
$\G$, and set $D := \mbox{max}\{\partial(x,y) \;|\; x,y \in X\}$. We call $D$  the {\it diameter} of 
$\G$. For a vertex $x \in X$ and an integer $i$ let $\G_i(x)$ denote the set of vertices at 
distance $i$ from $x$. We abbreviate $\G(x) = \G_1(x)$. For an integer $k \ge 0$ we say $\Gamma$ is 
{\em regular with valency} $k$ whenever $|\G(x)|=k$ for all $x \in X$. We say $\G$ is 
{\it distance-regular} whenever for all integers $h,i,j\;(0 \le h,i,j \le D)$ and for all vertices 
$x,y \in X$ with $\partial(x,y)=h,$ the number
\begin{eqnarray*}
p_{ij}^h = | \G_i(x) \cap \G_j(y)|
\end{eqnarray*}
is independent of $x$ and $y$. The $p_{ij}^h$ are called the {\it intersection numbers} of $\G$. 

\medskip \noindent
For the rest of this paper we assume $\Gamma$ is distance-regular with diameter $D \ge 3$. 
Note that $p^h_{ij}=p^h_{ji}$ for $0 \le h,i,j \le D$. For convenience set 
$c_i:=p_{1, i-1}^i \, (1 \le i \le D)$, $a_i:=p_{1i}^i \, (0 \le i \le D)$, 
$b_i:=p_{1, i+1}^i \, (0 \le i \le D-1)$, $k_i:=p_{ii}^0 \, (0 \le i \le D)$, and $c_0=b_D=0$. By the 
triangle inequality the following hold for $0 \le h,i,j \le D$: (i) $p_{ij}^h=0$ if one of $h,i,j$ is 
greater than the sum of the other two; (ii) $p_{ij}^h \ne 0$ if one of $h,i,j$ equals the sum of the 
other two. In particular $c_i \ne 0$ for $1 \le i \le D$ and $b_i \ne 0$ for $0 \le i \le D-1$. We 
observe that $\G$ is regular with valency $k=k_1=b_0$ and that 
\begin{equation}
\label{kk}
c_i+a_i+b_i=k \qquad \qquad (0 \le i \le D).
\end{equation}
Note that $k_i = |\G_i(x)|$ for $x \in X$ and $0 \le i \le D$. By \cite[p.~127]{BCN},
\begin{equation}
\label{k}
  k_i = \frac{b_0 b_1 \cdots b_{i-1}}{c_1 c_2 \cdots c_i} \qquad \qquad (0 \le i \le D).
\end{equation}

\medskip \noindent 
We recall the Bose-Mesner algebra of $\G$. For $0 \le i \le D$ let $A_i$ denote the matrix in 
$\MX$ with $(x,y)$-entry
\begin{equation}
\label{dm1}
  (A_i)_{x y} = \left\{ \begin{array}{ll}
                 1 & \hbox{if } \; \partial(x,y)=i, \\
                 0 & \hbox{if } \; \partial(x,y) \ne i \end{array} \right. \qquad (x,y \in X).
\end{equation}
We call $A_i$ the $i$th {\it distance matrix} of $\G$. We abbreviate $A:=A_1$ and call this the 
{\it adjacency matrix} of $\G.$ We observe
(ai)   $A_0 = I$;
(aii)  $\sum_{i=0}^D A_i = J$;
(aiii) $\overline{A_i} = A_i \;(0 \le i \le D)$;
(aiv)  $A_i^t = A_i  \;(0 \le i \le D)$;
(av)   $A_iA_j = \sum_{h=0}^D p_{ij}^h A_h \;(0 \le i,j \le D)$,
where $I$ (resp. $J$) denotes the identity matrix (resp. all 1's matrix) in  $\MX$. Using these facts
we find $A_0,A_1,\ldots,A_D$ is a basis for a commutative subalgebra $M$ of $\MX$. We call $M$ the 
{\it Bose-Mesner algebra} of $\G$. It turns out that $A$ generates $M$ \cite[p.~190]{BI}. By 
\cite[p.~45]{BCN}, $M$ has a second basis $E_0,E_1,\ldots,E_D$ such that 
(ei)   $E_0 = |X|^{-1}J$;
(eii)  $\sum_{i=0}^D E_i = I$;
(eiii) $\overline{E_i} = E_i \;(0 \le i \le D)$;
(eiv)  $E_i^t =E_i  \;(0 \le i \le D)$;
(ev)   $E_iE_j =\delta_{ij}E_i  \;(0 \le i,j \le D)$.
We call $E_0, E_1, \ldots, E_D $  the {\it primitive idempotents} of $\Gamma$.  

\medskip \noindent 
We now recall the Krein parameters. Let $\circ $ denote the entrywise product in $\MX$. Observe
$A_i\circ A_j= \delta_{ij}A_i$ for $0 \leq i,j\leq D$, so $M$ is closed under $\circ$. Thus there 
exist complex scalars $q^h_{ij}$  $(0 \leq h,i,j\leq D)$ such that
$$
  E_i \circ E_j = |X|^{-1}\sum_{h=0}^D q^h_{ij}E_h \qquad (0 \leq i,j\leq D).
$$
By \cite[Proposition 4.1.5]{BCN}, $q^h_{ij}$ is real and nonnegative  for $0 \leq h,i,j\leq D$. The 
$q^h_{ij}$ are called the {\it Krein parameters of} $\G$. The graph $\G$ is said to be 
{\it $Q$-polynomial} (with respect to the given ordering $E_0, E_1, \ldots, E_D$ of the primitive 
idempotents) whenever for $0 \leq h,i,j\leq D$, $q^h_{ij}= 0$ (resp. $q^h_{ij}\not= 0$) whenever one 
of $h,i,j$ is greater than (resp. equal to) the sum of the other two. 
For the rest of this section assume $\Gamma$ is $Q$-polynomial with respect to $E_0, E_1, \ldots, E_D$.

We now recall the dual idempotents of $\G$. To do this
fix a vertex $x \in X.$ We view $x$ as a ``base vertex''. For $ 0 \le i \le D$ let 
$E_i^*=E_i^*(x)$ denote the diagonal matrix in $\MX$ with $(y,y)$-entry 
\begin{equation}
\label{den0}
(\Es_i)_{y y} = \left\{ \begin{array}{lll}
                 1 & \hbox{if } \; \partial(x,y)=i, \\
                 0 & \hbox{if } \; \partial(x,y) \ne i \end{array} \right. 
                 \qquad (y \in X).
\end{equation}
We call $E_i^*$ the  $i$th 
{\it dual idempotent} of $\G$ with respect to $x$ \cite[p.~378]{ter1}. We observe
(i)   $\sum_{i=0}^D E_i^*=I$;
(ii)  $\overline{E_i^*} = E_i^*$ $(0 \le i \le D)$;
(iii) $E_i^{*t} = E_i^*$ $(0 \le i \le D)$;
(iv)  $E_i^*E_j^* = \delta_{ij}E_i^* $ $(0 \le i,j \le D)$.
By these facts $E_0^*,E_1^*, \ldots, E_D^*$ form a basis for a commutative subalgebra $M^*=M^*(x)$ 
of $\MX$. We call $M^*$ the {\it dual Bose-Mesner algebra} of $\G$ with respect to 
$x$ \cite[p.~378]{ter1}. For $0 \le i \le D$ we have
$$
\Es_i V = {\rm Span} \{ \h{y} \mid y \in X, \: \partial(x,y)=i\}
$$
so ${\rm dim} \Es_i V = k_i$. We call $\Es_i V$ the $i$th {\em subconstituent of} $\G$ 
with respect to $x$. Note that
\begin{equation}
\label{vsub}
V = E_0^*V + E_1^*V + \cdots + E_D^*V \qquad \qquad {\rm (orthogonal\ direct\ sum}).
\end{equation}
Moreover $\Es_i$ is the projection from $V$ onto $\Es_i V$ for $0 \le i \le D$.
Let $A^*=A^*(x)$ denote the diagonal matrix in $\MX$ with $(y,y)$-entry
$$
  A^*_{yy} = |X| E_{xy} \qquad (y \in X),
$$
where $E=E_1$. We call $A^*$ the {\em dual adjacency matrix of} $\G$ with respect to $x$. 
By \cite[Lemma 3.11(ii)]{ter1} $A^*$ generates $M^*$.

\medskip \noindent
We recall the Terwilliger algebra of $\G$. Let $T=T(x)$ denote the subalgebra of $\MX$ generated by 
$M$, $M^*$. We call $T$ the {\it Terwilliger algebra of} $\G$ with respect to $x$ 
\cite[Definition 3.3]{ter1}. Recall $M$ (resp. $M^*$) is generated by $A$ (resp. $\As$) so $T$ is 
generated by $A, A^*$. We observe $T$ has finite dimension. By construction $T$ is closed under the 
conjugate-transpose map so $T$ is semi-simple \cite[Lemma 3.4(i)]{ter1}. 

\medskip \noindent
By a $T$-{\em module} we mean a subspace $W$ of $V$ such that $BW \subseteq W$ 
for all $B \in T$. Let $W$ denote a $T$-module. Then $W$ is said to be {\em irreducible} whenever
$W$ is nonzero and $W$ contains no $T$-modules other than $0$ and $W$. Assume $W$ is irreducible.
Then $A$ and $\As$ act on $W$ as a tridiagonal pair \cite[Example 1.4]{ITT}.
We refer the reader to \cite{ITT, IT1, IT2, IT3, No, NT} and the references therein
for background on tridiagonal pairs.

By \cite[Corollary 6.2]{Go1} any $T$-module is an orthogonal direct sum of irreducible $T$-modules. 
In particular the standard module $V$ is an orthogonal direct sum of irreducible $T$-modules. Let 
$W, W'$ denote $T$-modules. By an {\em isomorphism of} $T$-{\em modules} from $W$ to $W'$ we mean an 
isomorphism of vector spaces $\sigma: W \to W'$ such that $(\sigma B - B \sigma) W = 0$ for all 
$B \in T$. The $T$-modules $W$, $W'$ are said to be {\em isomorphic} whenever there exists an 
isomorphism of $T$-modules from $W$ to $W'$. By \cite[Lemma 3.3]{Cu1} any two nonisomorphic 
irreducible $T$-modules are orthogonal. Let $W$ denote an irreducible $T$-module. By 
\cite[Lemma 3.4(iii)]{ter1} $W$ is an orthogonal direct sum of the nonvanishing spaces among 
$\Es_0 W, \Es_1 W, \ldots , \Es_D W$. By the {\em endpoint} of $W$ we mean 
$\min \{i \mid 0 \le i\le D, \; \Es_i W \ne 0 \}$. By the {\em diameter} of $W$ we mean 
$|\{i \mid 0 \le i \le D, \; \Es_i W \ne 0\}| -1$. 

By \cite[Proposition 8.3, Proposition 8.4]{egge1} $M \h{x}$ is the unique irreducible $T$-module with 
endpoint 0 and the unique irreducible $T$-module with diameter $D$. Moreover $M \h{x}$ is the unique 
irreducible $T$-module on which $E_0$ does not vanish. We call $M \h{x}$ the 
{\em primary module}. 

\medskip \noindent
We finish this section with some comments on local eigenvalues.
Let $\Delta=\Delta(x)$ denote the vertex-subgraph of $\Gamma$
induced on the set of vertices in $X$ adjacent $x$, and let $\breve{A}$ denote the adjacency matrix
of $\Delta$. By the {\em local eigenvalues of} $\G$ we mean the eigenvalues of $\breve{A}$. Note that
the local eigenvalues of $\G$ are precisely the eigenvalues of $\Es_1 A \Es_1$ on $\Es_1 V$.

Let $W$ denote an irreducible $T$-module with endpoint $1$.
By \cite[Theorem 2.2]{HI} $\Es_1 W$ is a one-dimensional eigenspace 
for $\Es_1 A \Es_1$; we call the corresponding 
eigenvalue the {\em local eigenvalue} of $W$.

\section{Distance-regular graphs of negative type}
\label{negtype}

In this section we recall what it means for
$\G$ to have classical parameters and negative type.
The graph $\G$ is said to have {\it classical parameters} $(D,b,\alpha,\beta)$ 
whenever the intersection numbers of $\G$ satisfy
$$
  c_i = {i \brack 1} \Big( 1+\alpha {i-1 \brack 1} \Big) \qquad \qquad (1 \le i \le D), 
$$
$$
  b_i = \Big( {D \brack 1} - {i \brack 1} \Big) 
        \Big( \beta-\alpha {i \brack 1} \Big) \qquad \qquad (0 \le i \le D-1), 
$$
where
$$
  {j \brack 1} := 1+b+b^2+ \cdots + b^{j-1}. 
$$
In this case $b$ is an integer and $b \not \in \{0,-1\}$.
If $\G$ has classical parameters then $\G$ is $Q$-polynomial \cite[Corollary 8.4.2]{BCN}.
We say that $\G$ has {\it negative type} whenever $\Gamma$ has classical parameters 
$(D,b,\alpha,\beta)$ such that $b < -1$. 

\medskip \noindent
We now recall kites and parallelograms.
Fix an integer $i \; (2 \le i \le D)$. By a {\it kite of length} $i$ (or 
$i$-{\it kite}) in $\G$ we mean a 4-tuple $uvwz$ of vertices of $\G$ 
such that $u, \, v, \, w$ are mutually adjacent, and $\partial(u,z) = i$, 
$\partial(v,z) = \partial(w,z)=i-1$. 
By a {\it parallelogram of length} $i$
(or $i$-{\it parallelogram}) in $\G$ we mean a 4-tuple $uvwz$ of vertices 
of $\G$ such that $\partial(u,v) = \partial(w,z)=1$, 
$\partial(u,z)=i,$ and $\partial(v,z)=\partial(u,w)=\partial(v,w)=i-1$. By \cite[Theorem 2.12]{ter2}
and \cite[Theorem 4.2]{Mi2}, if $\G$ has negative type then $\G$ has no parallelograms or
kites of any length.

\medskip \noindent
We now recall the near polygons. The graph $\G$ is called a {\em near polygon} whenever
$a_i=a_1 c_i$ for $1 \le i \le D-1$ and $\G$ has no $2$-kite \cite{No}. 
From now on we adopt the following notational convention.

\begin{notation}
\label{blank}
{\rm Let $\G=(X,R)$ denote a distance-regular graph with classical parameters $(D,b,\alpha,\beta)$,
$D \ge 3$, with valency $k$ and $a_1 \ne 0$. Assume that $\G$ is of negative type 
and $\G$ is not a near polygon. Let $A_0, A_1, \ldots, A_D$ denote the distance matrices of $\G$, 
and let $V$ denote the standard module of $\G$. We fix $x \in X$ and let 
$\Es_i=\Es_i(x) \; (0 \le i \le D)$, and $T=T(x)$ denote the corresponding dual idempotents and 
Terwilliger algebra, respectively.}
\end{notation}

\medskip \noindent
The following result is an immediate consequence of \cite[Lemma 6.4, Lemma 6.5]{Mi2}.

\begin{corollary}
\label{near}
With reference to Notation \ref{blank} we have $a_i > a_1 c_i$ for $2 \le i \le D$.
\end{corollary}

\section{The sets $\boldsymbol{D_j^i}$}
\label{sec:sets}

With reference to Notation \ref{blank}, in this section we define certain subsets $D^i_{\! j}$ of $X$ 
and explore their properties.

\begin{definition}
\label{def:D}
{\rm With reference to Notation \ref{blank} fix $z \in \G(x)$. For all integers $i,j$ 
we define $D^i_{\! j}=D_{\! j}^i(x,z)$ by
$$
  D^i_{\! j} = \G_i(x) \cap \G_j(z).
$$
We observe $D^i_{\! j} = \emptyset$ unless $0 \le i,j \le D$.
}
\end{definition} 

\begin{lemma}
\label{lem:D}
With reference to Notation \ref{blank} and Definition \ref{def:D} the following {\rm (i), (ii)} 
hold for $0 \le i,j \le D$.
\begin{itemize}
\item[{\rm (i)}] 
$|D^i_{\! j}|=p^1_{ij}$.
\item[{\rm (ii)}]
$D^i_{\! j} = \emptyset$ if and only if $p^1_{ij} = 0$.
\end{itemize}
\end{lemma}
\proof
(i) Immediate from the definition of $p^1_{ij}$ and $D^i_{\! j}$.

\smallskip \noindent
(ii) Immediate from (i) above. \qed

\begin{lemma}
\label{ps}
{\rm (\cite[p.~134]{BCN})}
With reference to Notation \ref{blank} the following {\rm (i), (ii)} hold.
\begin{itemize}
\item[{\rm (i)}] 
$p^1_{i-1,i} = p^1_{i,i-1} = c_i k_i k^{-1} \; (1 \le i \le D)$.
\item[{\rm (ii)}]
$p^1_{ii} = a_i k_i k^{-1} \; (0 \le i \le D)$.
\end{itemize}
\end{lemma}

\begin{lemma}
\label{ps1}
With reference to Notation \ref{blank} the following {\rm (i)--(iii)} hold.
\begin{itemize}
\item[{\rm (i)}] 
$p^1_{i-1,i} \ne 0$,  $p^1_{i,i-1} \ne 0 \; (1 \le i \le D)$.
\item[{\rm (ii)}]
$p^1_{00}=0$, $p^1_{ii} \ne 0 \; (1 \le i \le D)$.
\item[{\rm (iii)}]
$p^1_{ij} = 0$ if $|i-j| \not \in \{0,1\} \; (0 \le i,j \le D)$.
\end{itemize}
\end{lemma}
\proof
(i) Immediate from Lemma \ref{ps}(i).

\smallskip \noindent
(ii) It is clear that $p^1_{00}=0$ and $p^1_{11} = a_1 \ne 0$. Assume $2 \le i \le D$.
By Corollary \ref{near} we find $a_i \ne 0$ so $p^1_{ii} \ne 0$ in view of 
Lemma \ref{ps}(ii).
 
\smallskip \noindent
(iii) Immediate from the triangle inequality. \qed

\begin{lemma} 
\label{razdalje}
With reference to Notation \ref{blank} and Definition \ref{def:D} the following {\rm (i)--(iii)} hold.
\begin{itemize}
\item[{\rm (i)}]  
$\partial(u,y)=1$ for all distinct $u,y \in D_1^1$.
\item[{\rm (ii)}] 
There are no edges between $D_i^{i-1} \cup D_{i-1}^i$ and $D_{i-1}^{i-1}$ for $2 \le i \le D$.  
\item[{\rm (iii)}]
For $1 \le i \le D$ we have $\partial(u,y)=i$ for all $u \in D_1^1$ and all 
$y \in D_i^{i-1} \cup D_{i-1}^i$.
\end{itemize} 
\end{lemma}
\proof
(i) If $u, y$ are not adjacent, then $yxzu$ is a $2$-kite, a contradiction.

\smallskip \noindent
(ii) There does not exist adjacent vertices $v,w$ with $v \in D_{i-1}^i$ and  
$w \in D_{i-1}^{i-1}$; otherwise $xzwv$ is an $i$-parallelogram,
a contradiction. A similar argument shows that there 
does not exist adjacent vertices $v,w$ with $v \in D^{i-1}_i$ and  
$w \in D_{i-1}^{i-1}$. 

\smallskip \noindent
(iii) If $i = 1$ then the result is clear. Assume $2 \le i \le D$. By the triangle inequality we find
$\partial(u,y) \in \{i-1,i\}$. But $\partial(u,y) \ge i$ by (ii) above, 
and the result follows. \qed

\medskip \noindent
We end this section with two remarks on the local eigenvalues.

\begin{corollary} 
\label{local}
With reference to Notation \ref{blank} let $\Delta=\Delta(x)$ denote the vertex-subgraph of $\Gamma$
induced on the set of vertices in $X$ adjacent $x$. Then the following {\rm (i), (ii)} hold.
\begin{itemize}
\item[{\rm (i)}]  
$\Delta$ is a disjoint union of $k (a_1+1)^{-1}$ cliques, each consisting of $a_1+1$ vertices.
\item[{\rm (ii)}] 
The local eigenvalues of $\G$ are $a_1$ with multiplicity $k (a_1+1)^{-1}$, and $-1$
with multiplicity $k a_1 (a_1+1)^{-1}$.
\end{itemize}
\end{corollary}
\proof 
(i) Immediate from Lemma \ref{razdalje}(i),(ii). 

\smallskip \noindent
(ii) Immediate from (i) above. \qed

\begin{corollary} 
\label{local1}
With reference to Notation \ref{blank} let $W$ denote an irreducible $T$-module with endpoint $1$. 
Then the local eigenvalue of $W$ is $a_1$ or $-1$.
\end{corollary}
\proof
The local eigenvalue of $W$ is a local eigenvalue of $\G$.  
The result now follows from Corollary \ref{local}(ii). \qed

\section{The sets $\boldsymbol{{D_i^i}(0)}$ and $\boldsymbol{{D_i^i}(1)}$}
\label{sec:negtype}

With reference to Notation \ref{blank}, in this section we define certain subsets $D^i_i(0)$ and
$D_i^i(1)$ of $X$ and explore their properties.

\begin{lemma}
\label{lem:pom}
With reference to Notation \ref{blank} and Definition \ref{def:D}, for $1 \le i \le D$ and 
$y \in D_i^i$ we have $|\G_{i-1}(y) \cap D_1^1| \le 1$.
\end{lemma}
\proof Assume that $|\G_{i-1}(y) \cap D_1^1| \ge 2$ and pick distinct $u,v \in \G_{i-1}(y) \cap D_1^1$.
Then $xuvy$ is an $i$-kite, a contradiction. \qed

\begin{definition} 
\label{D()}
\rm 
With reference to Notation \ref{blank} and Definition \ref{def:D}, for an integer $i$ and 
$j \in \{0,1\}$ define a set $D_i^i(j)=D_i^i(j)(x,z)$ by
$$ 
  D_i^i(j) = \{y \in D_i^i \, | \, |\Gamma_{i-1}(y) \cap D_1^1| = j\}.
$$
We observe $D_i^i(j) = \emptyset$ unless $1 \le i \le D$.
By Lemma \ref{lem:pom} $D_i^i$ is the disjoint union of $D_i^i(1)$ and $D_i^i(0)$.
\end{definition}

\begin{lemma} 
\label{moc_1}
{\rm (\cite[Lemma 6.4]{Mi2})}
With reference to Notation \ref{blank} and Definition \ref{D()}
the following {\rm (i), (ii)} hold for $1 \le i \le D$.
\begin{itemize}
\item[{\rm (i)}]  $|D_i^i(1)|=a_1 c_i k_i k^{-1}$.
\item[{\rm (ii)}] $|D_i^i(0)|=(a_i-a_1 c_i) k_i k^{-1}$.
\end{itemize}
\end{lemma}

\begin{lemma}
\label{moc_2}
With reference to Notation \ref{blank} and Definition \ref{D()}
the following {\rm (i), (ii)} hold for $1 \le i \le D$.
\begin{itemize}
\item[{\rm (i)}]  $D_i^i(1) \ne \emptyset$ for $1 \le i \le D$.
\item[{\rm (ii)}] $D_i^i(0) \ne \emptyset$ for $2 \le i \le D$ and $D_1^1(0) = \emptyset$.
\end{itemize}
\end{lemma}
\proof Combine Corollary \ref{near} and Lemma \ref{moc_1}. \qed 

\begin{lemma} 
\label{razdalje2}
{\rm (\cite[Lemma 4.4(i)]{Mi2})}
With reference to Notation \ref{blank}, Definition \ref{def:D} and Definition \ref{D()},
for $2 \le i \le D$ we have $\partial(u,y)=i$ 
for all $u \in D_1^1$ and all $y \in D_i^i(0)$.             
\end{lemma}

\begin{lemma} \label{povezave}
{\rm (\cite[Sections 5,6]{Mi2})}
With reference to Notation \ref{blank}, Definition \ref{def:D} and Definition \ref{D()}
the following {\rm (i)--(iii)} hold.
\begin{itemize}
\item[{\rm (i)}] For  $1 \le i \le D$, each vertex in $D_{i-1}^i$ (resp. $D_i^{i-1}$) is adjacent to \\
    \begin{tabular}{ccl}     
     precisely & $c_{i-1}$ & vertices in $D_{i-2}^{i-1}$ 
                                 (resp. $D_{i-1}^{i-2}$), \\
     precisely & $c_i-c_{i-1}$ & vertices in $D_i^{i-1}$ 
                                     (resp. $D_{i-1}^i$), \\
     precisely & $a_{i-1}$ & vertices in $D_{i-1}^i$ 
                                 (resp. $D_i^{i-1}$), \\
     precisely & $b_i$ & vertices in $D_i^{i+1}$ (resp. $D_{i+1}^i$), \\
     precisely & $a_1(c_i-c_{i-1})$ & vertices in $D_i^i(1)$, \\
     precisely & $a_i - a_{i-1} - a_1(c_i-c_{i-1})$ & 
         vertices in $D_i^i(0)$,
    \end{tabular} \\
    and no other vertices in $X$.
\item[{\rm (ii)}] For  $2 \le i \le D$, each vertex in $D_i^i(0)$ is adjacent to \\
    \begin{tabular}{ccl}
     precisely & $c_i (b^{i-2}-1) (b^i-1)^{-1}$ & vertices in $D_{i-1}^{i-1}(0)$, \\
     precisely & $a_1 c_i (b^i - b^{i-2}) (b^i-1)^{-1}$ & vertices in $D_i^i(1)$, \\
     precisely & $c_i (b^i - b^{i-2}) (b^i-1)^{-1}$ & 
         vertices in $D_{i-1}^i$, \\
     precisely & $c_i (b^i - b^{i-2})(b^i-1)^{-1}$ &
         vertices in $D_i^{i-1}$, \\
     precisely & $b_i$ & vertices in $D_{i+1}^{i+1}(0)$, \\
     precisely & $a_i - c_i (a_1 + 1)(b^i - b^{i-2}) (b^i-1)^{-1}$ & 
         vertices in $D_i^i(0)$,
    \end{tabular} \\
    and no other vertices in $X$.
\item[{\rm (iii)}] For  $1 \le i \le D$, each vertex in $D_i^i(1)$ is adjacent to \\
    \begin{tabular}{ccl}
     precisely & $c_{i-1}$ & vertices in $D_{i-1}^{i-1}(1)$, \\
     precisely & $(a_1-1)(c_i-c_{i-1})+a_{i-1}$ & vertices in $D_i^i(1)$, \\
     precisely & $c_i-c_{i-1}$ & vertices in $D_{i-1}^i$, \\
     precisely & $c_i-c_{i-1}$ & vertices in $D_i^{i-1}$, \\
     precisely & $b_i$ & vertices in $D_{i+1}^{i+1}(1)$, \\
     precisely & $a_i - a_{i-1} - a_1(c_i-c_{i-1})$ & vertices in $D_i^i(0)$,
    \end{tabular} \\
    and no other vertices in $X$.
\end{itemize}
\end{lemma}

\section{Some products in $\boldsymbol{T}$}
\label{sec:prod}

With reference to Notation \ref{blank}, in this section we evaluate several products in $T$ which 
we shall need later. 

\begin{lemma}
\label{prod:lem0}
With reference to Notation \ref{blank}, for $0 \le h,i,j \le D$ and $y,z \in X$ the $(y,z)$-entry
of $\Es_h A_i \Es_j$ is $1$ if $\partial(x,y)=h$, $\partial(y,z)=i$, $\partial(x,z)=j$, and $0$
otherwise.
\end{lemma}
\proof
Compute the $(y,z)$-entry of $\Es_h A_i \Es_j$  by matrix multiplication and simplify the
result using \eqref{dm1} and \eqref{den0}. \qed

\begin{corollary}
\label{ter:lem1i}
{\rm (\cite[Lemma 3.2]{ter1})}
With reference to Notation \ref{blank},
$$
  \Es_h A_i \Es_j = 0 \quad \hbox{\rm if and only if } \quad p_{ij}^h=0 
  \qquad \qquad (0 \le h,i,j \le D).
$$
\end{corollary}
\proof
Immediate from Lemma \ref{prod:lem0}. \qed

\begin{corollary}
\label{ter:cor1i}
With reference to Notation \ref{blank} and Definition \ref{def:D},
for $0 \le i,j \le D$ and $y \in X$ the $(y,z)$-entry
of $\Es_i A_j \Es_1$ is $1$ if $y \in D^i_{\! j}$, and $0$ otherwise.
\end{corollary}
\proof
Immediate from Lemma \ref{prod:lem0} . \qed

\begin{corollary}
\label{prod:cor0}
With reference to Notation \ref{blank} the following {\rm (i), (ii)} hold.
\begin{itemize}
\item[{\rm (i)}]
$\Es_i A_{i-1} \Es_1 + \Es_i A_i \Es_1 + \Es_i A_{i+1} \Es_1 = \Es_i J \Es_1$ for $1 \le i \le D-1$.
\item[{\rm (ii)}]
$\Es_D A_{D-1} \Es_1 + \Es_D A_D \Es_1 = \Es_D J \Es_1$.
\end{itemize}
\end{corollary}
\proof
For each equation evaluate the right-hand side using assertion (aii) below line \eqref{dm1}, 
and simplify the result using Corollary \ref{ter:lem1i} and assertion (i) above line \eqref{kk}. \qed

\begin{lemma}
\label{prod:lem00}
With reference to Notation \ref{blank}, for $0 \le h,i,j,r,s \le D$ and $y,z \in X$ the $(y,z)$-entry
of $\Es_h A_r \Es_i A_s \Es_j$ is $|\G_i(x) \cap \G_r(y) \cap \G_s(z)|$ if $\partial(x,y)=h$,
$\partial(x,z)=j$, and $0$ otherwise.
\end{lemma}
\proof
Compute the $(y,z)$-entry of $\Es_h A_r \Es_i A_s \Es_j$  by matrix multiplication and simplify the
result using \eqref{dm1} and \eqref{den0}. \qed

\begin{corollary}
\label{prod:cor1}
With reference to Notation \ref{blank} and Definition \ref{D()},
for $1 \le i \le D$ and $y \in X$ the following {\rm (i), (ii)} hold.
\begin{itemize}
\item[{\rm (i)}] The $(y,z)$-entry of $\Es_i A_{i-1} \Es_1 A \Es_1$ is $1$ if $y \in D_i^i(1)$, 
and $0$ otherwise.
\item[{\rm (ii)}] The $(y,z)$-entry of $\Es_i (A_i - A_{i-1} \Es_1 A) \Es_1$ is $1$ if $y \in D_i^i(0)$, 
and $0$ otherwise.
\end{itemize}
\end{corollary}
\proof
(i) Immediate from Lemma \ref{razdalje}(iii), Lemma \ref{razdalje2} and Lemma \ref{prod:lem00}. 

\smallskip \noindent
(ii) By Corollary \ref{ter:cor1i} the $(y,z)$-entry of $\Es_i A_i \Es_1$ is $1$ if 
$y \in D_i^i$, and $0$ otherwise. The result now follows from (i) above
and since $D_i^i$ is the disjoint union of $D_i^i(0)$ and $D_i^i(1)$. \qed

\section{The matrices $\boldsymbol{L}$, $\boldsymbol{F}$, $\boldsymbol{R}$}
\label{sec:LFR}

With reference to Notation \ref{blank}, in this section we recall the matrices $L$, $F$, $R$
and use them to interpret Theorem \ref{povezave}.
 
\begin{definition}
\label{def:LFR}
{\rm With reference to Notation \ref{blank} we define matrices $L=L(x)$, $F=F(x)$, $R=R(x)$ by
$$
  L = \sum_{h=1}^D \Es_{h-1} A \Es_h, \qquad \qquad
  F = \sum_{h=0}^D \Es_h A \Es_h, \qquad \qquad
  R = \sum_{h=0}^{D-1} \Es_{h+1} A \Es_h.
$$
Note that $A=L+F+R$ \cite[Lemma 4.4]{Cu1}. We call $L$, $F$, and $R$ the {\em lowering matrix}, 
the {\em flat matrix}, and the {\em raising matrix of} $\G$ with respect to $x$.}
\end{definition}

\begin{lemma}
\label{prod:lem1}
With reference to Notation \ref{blank} and Definition \ref{def:LFR}
the following {\rm (i)--(iii)} hold.
\begin{itemize}
\item[{\rm (i)}] 
$L \Es_1 = \Es_0 A \Es_1$.
\item[{\rm (ii)}] 
For $2 \le i \le D$, 
$$
  L \Es_i A_{i-1} \Es_1 = b_{i-1} \Es_{i-1} A_{i-2} \Es_1 + 
                        (c_i - c_{i-1}) \Es_{i-1} A_i \Es_1.
$$                       
\item[{\rm (iii)}]
For $1 \le i \le D-1$,
$$
  L \Es_i A_{i+1} \Es_1 = b_i \Es_{i-1} A_i \Es_1.
$$
\end{itemize}
\end{lemma}
\proof 
For each equation and for $y,z \in X$ compute the $(y,z)$-entry of each side and interpret the results 
using Theorem \ref{povezave}, Corollary \ref{ter:cor1i} and Lemma \ref{prod:lem00}. \qed

\begin{lemma}
\label{prod:lem2}
With reference to Notation \ref{blank} and Definition \ref{def:LFR} 
the following {\rm (i)--(iii)} hold.
\begin{itemize}
\item[{\rm (i)}]
$F \Es_1 = \Es_1 A \Es_1$.
\item[{\rm (ii)}] 
For $2 \le i \le D$,
\begin{eqnarray*}
F \Es_i A_{i-1} \Es_1 & = & a_{i-1} \Es_i A_{i-1} \Es_1 + 
                            (c_i - c_{i-1}) \Es_i A_{i-1} \Es_1 A \Es_1 \\
                      &   & + \, c_i (b^i - b^{i-2}) (b^i-1)^{-1} \Es_i (A_i - A_{i-1} \Es_1 A) \Es_1.
\end{eqnarray*}
\item[{\rm (iii)}]
For $1 \le i \le D-1$,
$$
  F \Es_i A_{i+1} \Es_1 = a_i \Es_i A_{i+1} \Es_1.
$$
\end{itemize}
\end{lemma}
\proof 
For each equation and for $y,z \in X$ compute the $(y,z)$-entry of each side and interpret the results 
using Theorem \ref{povezave}, Corollary \ref{ter:cor1i}, Lemma \ref{prod:lem00} and Corollary \ref{prod:cor1}. \qed

\begin{lemma}
\label{prod:lem3}
With reference to Notation \ref{blank} and Definition \ref{def:LFR} 
the following {\rm (i)--(iv)} hold.
\begin{itemize}
\item[{\rm (i)}] 
For $1 \le i \le D-1$,
$$
  R \Es_i A_{i-1} \Es_1 = c_i \Es_{i+1} A_i \Es_1.
$$
\item[{\rm (ii)}]
$R \Es_D A_{D-1} \Es_1 = 0$.
\item[{\rm (iii)}]
For $1 \le i \le D-2$,
\begin{eqnarray*}
R \Es_i A_{i+1} \Es_1 & = & c_{i+1} \Es_{i+1} A_{i+2} \Es_1 + (c_{i+1}-c_i) \Es_{i+1} A_i \Es_1 \\
                      &   & + \, (c_{i+1}-c_i) \Es_{i+1} A_i \Es_1 A \Es_1 \\
                      &   & + \, c_{i+1} (b^{i+1} - b^{i-1}) (b^{i+1}-1)^{-1} \Es_{i+1}(A_{i+1}-A_i \Es_1 A)\Es_1.
\end{eqnarray*}                         
\item[{\rm (iv)}]
\begin{eqnarray*}
R \Es_{D-1} A_D \Es_1 & = & (c_D-c_{D-1}) \Es_D A_{D-1} \Es_1 + (c_D-c_{D-1}) \Es_D A_{D-1} \Es_1 A \Es_1 \\
                      &   & + \, c_D (b^D - b^{D-2}) (b^D-1)^{-1} \Es_D(A_D-A_{D-1} \Es_1 A)\Es_1.                    
\end{eqnarray*}                            
\end{itemize}
\end{lemma}
\proof 
For each equation and for $y,z \in X$ compute the $(y,z)$-entry of each side and interpret the results 
using Theorem \ref{povezave}, Corollary \ref{ter:cor1i}, Lemma \ref{prod:lem00} and Corollary \ref{prod:cor1}. \qed

\section{More products in $\boldsymbol{T}$}
\label{sec:more}

With reference to Notation \ref{blank}, in this section we evaluate more products in $T$ 
which we will need later.

\begin{lemma}
\label{prod:lem4A}
With reference to Notation \ref{blank}, for $y,z \in \G(x)$ and $1 \le i \le D$ the number
$|\G_i(x) \cap \G_{i-1}(y) \cap \G_{i-1}(z)|$ is equal to 
$c_i k_i k^{-1}$ if $y=z$, 
$0$ if $\partial(y,z)=1$, and 
$c_i (c_i-1) k_i k^{-1} b_1^{-1}$ if $\partial(y,z) = 2$.
\end{lemma}
\proof 
If $y=z$ then the result follows by Lemma \ref{ps}(i). 
If $\partial(y,z)=1$ then the result follows by Lemma \ref{razdalje}(iii). 
Assume $\partial(y,z)=2$. Abbreviate $D_{\! j}^{\ell}=D_{\! j}^{\ell}(x,z) \; (0 \le j, \ell \le D)$ 
and note that $y \in D_2^1$. It follows from Theorem \ref{povezave} that the number of paths of 
length $i-1$ between $y$ and $D^i_{i-1}$ is independent of $y$. Moreover, between any two 
vertices of $\G$ which are at distance $i-1$, there exist exactly $c_1 c_2 \cdots c_{i-1}$ paths 
of length $i-1$. Therefore the scalar $|D^i_{i-1} \cap \G_{i-1}(y)|$ is independent of $y$; denote
this scalar by $\alpha_i$. For $v \in D^i_{i-1}$ we have $|\G_{i-1} (v) \cap \G(x)| = c_i$, 
so using Lemma \ref{razdalje}(iii) we find $|\G_{i-1} (v) \cap D^1_2| = c_i-1$. 
Using these comments we count in two ways the number of pairs 
$(y,v)$ such that $y \in D^1_2$, $v \in D^i_{i-1}$, and $\partial(y,v)=i-1$.
This yields $\alpha_i |D^1_2| = |D^i_{i-1}|(c_i-1)$.
Evaluating this equation using Lemma \ref{lem:D}(i) and Lemma \ref{ps}(i) we find
$\alpha_i = c_i (c_i-1) k_i k^{-1} b_1^{-1}$.
The result follows. \qed

\begin{corollary}
\label{prod:corA}
With reference to Notation \ref{blank}, for $1 \le i \le D$ we have
$$
  \Es_1 A_{i-1} \Es_i A_{i-1} \Es_1 = c_i k_i k^{-1} \Es_1 + 
                                     c_i (c_i-1) k_i k^{-1} b_1^{-1} \Es_1 A_2 \Es_1.
$$
\end{corollary}
\proof For $y,z \in X$ we show that the $(y,z)$-entry of both sides are equal. 
If $y \not \in \G(x)$ or $z \not \in \G(x)$ then the $(y,z)$-entry of each side is $0$.
If $y,z \in \G(x)$ then the $(y,z)$-entry of both sides are equal by Corollary \ref{ter:cor1i},
Lemma \ref{prod:lem00} and Lemma \ref{prod:lem4A}. The result follows. \qed

\begin{lemma}
\label{prod:lem4F}
With reference to Notation \ref{blank}, for $y,z \in \G(x)$ and $1 \le i \le D-1$ the number
$|\G_i(x) \cap \G_{i-1}(y) \cap \G_{i+1}(z)|$ is equal to 
$0$ if $y=z$, 
$0$ if $\partial(y,z)=1$, and 
$c_i b_i k_i k^{-1} b_1^{-1}$ if $\partial(y,z)=2$.
\end{lemma}
\proof 
If $y=z$ then the result is clear. 
If $\partial(y,z)=1$ then the result follows by Lemma \ref{razdalje}(iii). 
Assume $\partial(y,z)=2$. Abbreviate $D_{\! j}^{\ell}=D_{\! j}^{\ell}(x,z) \; (0 \le j, \ell \le D)$ 
and note that $y \in D_2^1$.  It follows from Theorem \ref{povezave} that the number of paths 
of length $i-1$ between $y$ and $D^i_{i+1}$ is independent of $y$. 
Moreover, between any two vertices of $\Gamma$ which are at distance 
$i-1$, there exist exactly $c_1 c_2 \cdots c_{i-1}$ paths of length $i-1$.
Therefore the scalar $|D^i_{i+1} \cap \G_{i-1}(y)|$ is independent of $y$; denote this scalar
by $\alpha_i$. For $v \in D^i_{i+1}$ we have 
$|\G_{i-1} (v) \cap \G(x)| = c_i$, so using Lemma \ref{razdalje}(iii) we find
$|\G_{i-1} (v) \cap D^1_2| = c_i$. 
Using these comments we count in two ways the number of pairs 
$(y,v)$ such that $y \in D^1_2$, $v \in D^i_{i+1}$, and $\partial(y,v)=i-1$.
This yields $\alpha_i |D^1_2| = |D^i_{i+1}|c_i$.
Evaluating this equation using Lemma \ref{lem:D}(i), Lemma \ref{ps}(i) and $c_{i+1} k_{i+1} = b_i k_i$ we find
$\alpha_i = c_i b_i k_i k^{-1} b_1^{-1}$.
The result follows. \qed

\begin{corollary}
\label{prod:corF}
With reference to Notation \ref{blank}, for $1 \le i \le D-1$ we have
$$
  \Es_1 A_{i-1} \Es_i A_{i+1} \Es_1 = c_i b_i k_i k^{-1} b_1^{-1} \Es_1 A_2 \Es_1.
$$
\end{corollary}
\proof For $y,z \in X$ we show that the $(y,z)$-entry of both sides are equal. 
If $y \not \in \G(x)$ or $z \not \in \G(x)$ then the $(y,z)$-entry of each side is $0$.
If $y,z \in \G(x)$ then the $(y,z)$-entry of both sides are equal by Corollary \ref{ter:cor1i},
Lemma \ref{prod:lem00} and Lemma \ref{prod:lem4F}. The result follows. \qed

\begin{lemma}
\label{prod:lem4E}
With reference to Notation \ref{blank}, for $y,z \in \G(x)$ and $1 \le i \le D-1$ the number
$|\G_i(x) \cap \G_{i+1}(y) \cap \G_{i+1}(z)|$ is equal to 
$b_i k_i k^{-1}$ if $y=z$, 
$b_i k_i k^{-1}$ if $\partial(y,z)=1$, and 
$b_i (b_1-a_i-c_i) k_i k^{-1} b_1^{-1}$ if $\partial(y,z)=2$.
\end{lemma}
\proof 
If $y=z$ the result follows by Lemma \ref{ps}(i) and since $c_{i+1} k_{i+1} = b_i k_i$. 
If $\partial(y,z)=1$ then the result follows by Lemma \ref{ps}(i), Lemma \ref{razdalje}(iii) and since 
$c_{i+1} k_{i+1} = b_i k_i$. 
Assume $\partial(y,z)=2$. Abbreviate $D_{\! j}^{\ell}=D_{\! j}^{\ell}(x,z) \; (0 \le j, \ell \le D)$ 
and note that $y \in D_2^1$. We first claim that $|D^i_{i+1} \cap \G_i(y)| = a_i b_i k_i k^{-1} b_1^{-1}$.
It follows from Theorem \ref{povezave} that the number of paths 
of length $i$ between $y$ and $D^i_{i+1}$ is independent of $y$. 
Moreover, between any two vertices of $\Gamma$ which are at distance 
$i-1$ ($i$, respectively), there exist exactly $(a_1 + a_2 + \cdots + a_{i-1}) c_1 c_2 \cdots c_{i-1}$
($c_1 c_2 \cdots c_i$, respectively) paths of length $i$. By this and Lemma \ref{prod:lem4F}
the scalar $|D^i_{i+1} \cap \G_i(y)|$ is independent of $y$; denote this scalar by $\alpha_i$.
For $v \in D_{i+1}^i$ we have $|\G_i (v) \cap \G(x)| = a_i$, so using Lemma \ref{razdalje}(iii) we find
$|\G_i (v) \cap D^1_2| = a_i$. 
Using these comments we count in two ways the number of pairs 
$(y,v)$ such that $y \in D^1_2$, $v \in D^i_{i+1}$, and $\partial(y,v)=i$.
This yields $\alpha_i |D^1_2| = |D^i_{i+1}| \, a_i$.
Evaluating this equation using Lemma \ref{lem:D}(i), Lemma \ref{ps}(i) and $c_{i+1}k_{i+1}=b_i k_i$ we find
\begin{equation}
\label{lem4E:pom}
 \alpha_i = a_i b_i k_i k^{-1} b_1^{-1}.
\end{equation}
We have proved the claim. We can now easily show that 
$|D^i_{i+1} \cap \G_{i+1}(y)| = b_i (b_1-a_i-c_i) k_i k^{-1} b_1^{-1}$.
Pick $v \in D^i_{i+1}$. It follows from the triangle inequality that 
$\partial(y,v) \in \{i-1,i,i+1\}$, so 
$$
  |D^i_{i+1} \cap \G_{i+1}(y)| = |D_{i+1}^i| - |D^i_{i+1} \cap \G_i(y)| - |D^i_{i+1} \cap \G_{i-1}(y)|.
$$
Using Lemma \ref{lem:D}(i), Lemma \ref{ps}(i), Lemma \ref{prod:lem4F} and \eqref{lem4E:pom} we find
$|D^i_{i+1} \cap \G_{i+1}(y)| = b_i (b_1-a_i-c_i) k_i k^{-1} b_1^{-1}$. The result follows. \qed

\begin{corollary}
\label{prod:corE}
With reference to Notation \ref{blank}, for $1 \le i \le D-1$ we have
$$
  \Es_1 A_{i+1} \Es_i A_{i+1} \Es_1 = b_i k_i k^{-1} \Es_1 + b_i k_i k^{-1} \Es_1 A \Es_1 +
                                      b_i (b_1-a_i-c_i) k_i k^{-1} b_1^{-1} \Es_1 A_2 \Es_1.
$$
\end{corollary}
\proof For $y,z \in X$ we show that the $(y,z)$-entry of both sides are equal. 
If $y \not \in \G(x)$ or $z \not \in \G(x)$ then the $(y,z)$-entry of each side is $0$.
If $y,z \in \G(x)$ then the $(y,z)$-entry of both sides are equal by Corollary \ref{ter:cor1i},
Lemma \ref{prod:lem00} and Lemma \ref{prod:lem4E}. The result follows. \qed

\section{Some scalar products}
\label{sec:sp}

With reference to Notation \ref{blank}, in this section we compute some scalar products which 
we will need later.

\begin{lemma}
\label{sp:lem0}
With reference to Notation \ref{blank} let $W$ denote an irreducible $T$-module with endpoint $1$.
Then $JW=0$.
\end{lemma}
\proof Since $W$ is not the primary module we have $E_0 W =0$. Recall $J=|X| E_0$ so $J W = 0$. \qed

\begin{lemma}
\label{sp:lem1}
With reference to Notation \ref{blank} let $W$ denote an irreducible $T$-module with endpoint $1$.
Then the following {\rm (i), (ii)} hold for $w \in \Es_1 W$.
\begin{itemize}
\item[{\rm (i)}]
$\Es_i A_{i-1} w + \Es_i A_i w + \Es_i A_{i+1} w = 0$ for $1 \le i \le D-1$.
\item[{\rm (ii)}]
$\Es_D A_{D-1} w + \Es_D A_D w = 0$.
\end{itemize}
\end{lemma}
\proof
For each equation in Corollary \ref{prod:cor0} apply both sides to $w$ and simplify using
$\Es_1 w = w$ and Lemma \ref{sp:lem0}. \qed

\begin{corollary}
\label{sp:cor1a}
With reference to Notation \ref{blank} let $W$ denote an irreducible $T$-module with endpoint $1$
and local eigenvalue $\eta$. Then for $w \in \Es_1 W$ we have $\Es_1 A_2 w = -(1+\eta) w$.
\end{corollary}
\proof
Set $i=1$ in Lemma \ref{sp:lem1}(i) and note that $\Es_1 A_0 w = w$ and $\Es_1 A w = \eta w$. \qed

\begin{lemma}
\label{sp:lem2}
With reference to Notation \ref{blank} let $W$ denote an irreducible $T$-module with endpoint $1$
and local eigenvalue $\eta$. Then the following {\rm (i)--(iii)} hold for $w \in \Es_1 W$.
\begin{itemize}
\item[{\rm (i)}]
$\Vert \Es_i A_{i-1} w \Vert^2 = (b_1-(c_i-1)(1+ \eta)) c_i k_i k^{-1} b_1^{-1} \Vert w \Vert^2 \; \; (1 \le i \le D)$.
\item[{\rm (ii)}]
$\Vert \Es_i A_{i+1} w \Vert^2 = 
(k-b_i)(1+\eta) b_i k_i k^{-1} b_1^{-1} \Vert w \Vert^2 \; \; (1 \le i \le D-1)$.
\item[{\rm (iii)}]
$\la \Es_i A_{i-1} w, \Es_i A_{i+1} w \ra = 
- (1+ \eta) c_i b_i k_i k^{-1} b_1^{-1} \Vert w \Vert^2 \; \; (1 \le i \le D-1)$.
\end{itemize}
\end{lemma}
\proof 
(i) Evaluating $\Vert \Es_i A_{i-1} w \Vert^2 = \la  \Es_i A_{i-1} w,  \Es_i A_{i-1} w \ra$ 
using $\Es_1 w = w$, line \eqref{cez} and Corollary \ref{prod:corA} we find 
$$
  \Vert \Es_i A_{i-1} w \Vert^2 = 
  \frac{c_i k_i}{k} \Vert w \Vert^2 + \frac{c_i (c_i-1) k_i}{k b_1} \la w, \Es_1 A_2 w \ra.
$$
The result follows from this and Corollary \ref{sp:cor1a}.

\noindent
(ii),(iii) Similar to the proof of (i) above. \qed

\medskip \noindent
We now split the analysis into two cases, depending on whether $W$ has local eigenvalue $-1$ or $a_1$.

\section{The irreducible $\boldsymbol{T}$-modules with endpoint $\boldsymbol{1}$ 
         and local eigenvalue $\boldsymbol{-1}$}
\label{sec:jeme}

With reference to Notation \ref{blank}, in this section we describe the irreducible $T$-modules 
with endpoint $1$ and local eigenvalue $-1$. 

\begin{lemma}
\label{jeme:lem1}
With reference to Notation \ref{blank} let $W$ denote an irreducible $T$-module with endpoint $1$
and local eigenvalue $-1$. Then for $w \in \Es_1 W$ and $1 \le i \le D-1$ we have
$\Es_i A_{i+1} w = 0$.
\end{lemma}
\proof Immediate from Lemma \ref{sp:lem2}(ii). \qed

\begin{theorem}
\label{jeme:thm5}
With reference to Notation \ref{blank} let $W$ denote an irreducible $T$-module with endpoint $1$
and local eigenvalue $-1$. Fix a nonzero $w \in \Es_1 W$. Then the following is a basis for $W$:
\begin{equation}
\label{basisjeme}
  \Es_i A_{i-1} w \qquad \qquad (1 \le i \le D).
\end{equation}
\end{theorem}
\proof
We first show that $W$ is spanned by the vectors \eqref{basisjeme}. Let $W'$ denote the subspace of $V$ 
spanned by the vectors \eqref{basisjeme} and note that $W' \subseteq W$. We claim that $W'$ is a $T$-module. 
By construction $W'$ is $M^*$-invariant. It follows from Lemma \ref{prod:lem1}(i),(ii), Lemma \ref{prod:lem2}(i),(ii), 
Lemma \ref{prod:lem3}(i),(ii), Lemma \ref{sp:lem1}, Lemma \ref{jeme:lem1},
$\Es_1 w = w$ and $\Es_1 A w = - w$ that $W'$ is 
invariant under each of $L$, $F$, $R$. Recall that $L+F+R=A$ and $A$ generates $M$ so $W'$ is $M$-invariant. 
The claim follows. Note that $W' \ne 0$ since $w \in W'$ so $W'=W$ by the irreducibility of $W$. 
We now show that the vectors \eqref{basisjeme} are linearly independent. By \eqref{vsub}
it suffices to show that $\Es_i A_{i-1} w \ne 0$ for $1 \le i \le D$. This follows from Lemma 
\ref{sp:lem2}(i) and since $w \ne 0$. \qed

\begin{corollary}
\label{jeme:cor6}
With reference to Notation \ref{blank} let $W$ denote an irreducible $T$-module with endpoint $1$
and local eigenvalue $-1$. Then $\Es_i W$ has dimension $1$ for $1 \le i \le D$.
\end{corollary}
\proof Immediate from Theorem \ref{jeme:thm5}. \qed

\begin{corollary}
\label{jeme:cor7}
With reference to Notation \ref{blank} let $W$ denote an irreducible $T$-module with endpoint $1$
and local eigenvalue $-1$. Then the following {\rm (i), (ii)} hold.
\begin{itemize}
\item[{\rm (i)}]
The dimension of $W$ is $D$.
\item[{\rm (ii)}]
The diameter of $W$ is $D-1$.
\end{itemize}
\end{corollary}
\proof Immediate from Corollary \ref{jeme:cor6}. \qed

\begin{corollary}
\label{jeme:cor8}
With reference to Notation \ref{blank} let $W$ denote an irreducible $T$-module with endpoint $1$
and local eigenvalue $-1$. Then $W=M^*Mw$ for all nonzero $w \in \Es_1 W$.
\end{corollary}
\proof By construction $M^*M w \subseteq W$ and equality holds in view of Theorem \ref{jeme:thm5}. \qed


\section{The irreducible $\boldsymbol{T}$-modules with endpoint $\boldsymbol{1}$ 
         and local eigenvalue $\boldsymbol{-1}$: the $\boldsymbol{A}$-action}
\label{sec:jemeact}

With reference to Notation \ref{blank} let $W$ denote an irreducible $T$-module with endpoint 1
and local eigenvalue $-1$.
In this section we display the action of $A$ on the basis for $W$ given in Theorem \ref{jeme:thm5}.
Since $A=L+F+R$ it suffices to give the actions of $L$, $F$, $R$ on this basis.

\begin{lemma}
\label{jemeact:lem1}
With reference to Notation \ref{blank} let $W$ denote an irreducible $T$-module with endpoint $1$
and local eigenvalue $-1$. Then the following {\rm (i), (ii)} hold for all nonzero $w \in \Es_1 W$.
\begin{itemize}
\item[{\rm (i)}]
$L w = 0$.
\item[{\rm (ii)}]
For $2 \le i \le D$,
$$
  L \Es_i A_{i-1} w = b_{i-1} \Es_{i-1} A_{i-2} w.
$$
\end{itemize}
\end{lemma}
\proof
For each equation of Lemma \ref{prod:lem1}(i),(ii) apply each side to $w$ and simplify using $\Es_1 w = w$
and Lemma \ref{jeme:lem1}. \qed

\begin{lemma}
\label{jemeact:lem2}
With reference to Notation \ref{blank} let $W$ denote an irreducible $T$-module with endpoint $1$
and local eigenvalue $-1$. Then the following holds for all nonzero $w \in \Es_1 W$ and $1 \le i \le D$:
$$
  F \Es_i A_{i-1} w = (a_{i-1}+c_{i-1}-c_i) \Es_i A_{i-1} w.
$$
\end{lemma}
\proof
For each equation of Lemma \ref{prod:lem2}(i),(ii) apply each side to $w$ and simplify using $\Es_1 w = w$,
$\Es_1 A w = - w$, Lemma \ref{sp:lem1} and Lemma \ref{jeme:lem1}. \qed

\begin{lemma}
\label{jemeact:lem3}
With reference to Notation \ref{blank} let $W$ denote an irreducible $T$-module with endpoint $1$
and local eigenvalue $-1$. Then the following {\rm (i), (ii)} hold for all nonzero $w \in \Es_1 W$.
\begin{itemize}
\item[{\rm (i)}]
For $1 \le i \le D-1$,
$$
  R \Es_i A_{i-1} w = c_i \Es_{i+1} A_i w.
$$
\item[{\rm (ii)}]
$R \Es_D A_{D-1} w = 0$.
\end{itemize}
\end{lemma}
\proof
For each equation of Lemma \ref{prod:lem3}(i),(ii) apply each side to $w$ and simplify using 
$\Es_1 w = w$. \qed


\section{The irreducible $\boldsymbol{T}$-modules with endpoint $\boldsymbol{1}$ 
         and local eigenvalue $\boldsymbol{a_1}$}
\label{sec:nime}

With reference to Notation \ref{blank}, in this section we describe the irreducible $T$-modules 
with endpoint $1$ and local eigenvalue $a_1$. 

\begin{lemma}
\label{nime:lem1}
With reference to Notation \ref{blank} let $W$ denote an irreducible $T$-module with endpoint $1$
and local eigenvalue $a_1$. For $w \in \Es_1 W$ and $2 \le i \le D-1$ the determinant of 
$$
\begin{pmatrix}
  \Vert \Es_i A_{i-1} w \Vert^2 & \la \Es_i A_{i-1} w, \Es_i A_{i+1} w \ra \\ 
  \la \Es_i A_{i+1} w, \Es_i A_{i-1} w \ra & \Vert \Es_i A_{i+1} w \Vert^2
\end{pmatrix}
$$
is equal to
$$
  c_i b_i (a_1+1)(a_i-a_1 c_i) k_i^2 k^{-1} b_1^{-2} \Vert w \Vert^4.
$$
\end{lemma}
\proof 
Evaluate the matrix entries using Lemma \ref{sp:lem2}, take the determinant and
simplify the result using \eqref{kk}. \qed

\begin{theorem}
\label{nime:thm2}
With reference to Notation \ref{blank} let $W$ denote an irreducible $T$-module with endpoint $1$
and local eigenvalue $a_1$. Fix a nonzero $w \in \Es_1 W$. Then the following is a basis for $W$:
\begin{equation}
\label{basisnime}
  \Es_i A_{i-1} w \quad (1 \le i \le D), \qquad \qquad
  \Es_i A_{i+1} w \quad (2 \le i \le D-1).
\end{equation}
\end{theorem}
\proof
We first show that $W$ is spanned by the vectors \eqref{basisnime}. Let $W'$ denote the subspace of $V$ 
spanned by the vectors \eqref{basisnime} and note that $W' \subseteq W$. We claim that $W'$ is a 
$T$-module. By construction $W'$ is $M^*$-invariant. It follows from Lemmas \ref{prod:lem1}, 
\ref{prod:lem2}, \ref{prod:lem3}, Lemma \ref{sp:lem1}, Corollary \ref{sp:cor1a}, $\Es_1 w = w$ and
$\Es_1 A w = a_1 w$ that $W'$ is invariant under each of 
$L$, $F$, $R$. Recall that $L+F+R=A$ and $A$ generates $M$ so $W'$ is $M$-invariant. The claim 
follows. Note that $W' \ne 0$ since $w \in W'$ so $W'=W$ by the irreducibility of $W$. 

We now show that the vectors \eqref{basisnime} are linearly independent. 
By \eqref{vsub} and since $w \ne 0$, it suffices to show that
$\Es_i A_{i-1} w$, $\Es_i A_{i+1}$ 
are linearly independent for $2 \le i \le D-1$, and that $\Es_D A_{D-1} w \ne 0$. 
For  $2 \le i \le D-1$ the vectors 
$\Es_i A_{i-1} w$, $\Es_i A_{i+1} w$ are linearly independent since their matrix of inner products
has nonzero determinant by Corollary \ref{near} and Lemma \ref{nime:lem1}. 
It follows from \eqref{kk} and Lemma \ref{sp:lem2}(i) that 
$\Vert \Es_D A_{D-1} w \Vert^2 = (a_D - a_1 c_D) c_D k_D k^{-1} b_1^{-1} \Vert w \Vert^2$. 
Now $\Es_D A_{D-1} w \ne 0$ by Corollary \ref{near}. 
By these comments the vectors \eqref{basisnime} are linearly independent and the result follows. \qed 

\medskip \noindent
We emphasize an idea from the proof of Theorem \ref{nime:thm2}.

\begin{corollary}
\label{nime:cor3}
With reference to Notation \ref{blank} let $W$ denote an irreducible $T$-module with endpoint $1$
and local eigenvalue $a_1$. Fix a nonzero $w \in \Es_1 W$. Then the following {\rm (i)--(iii)} hold.
\begin{itemize}
\item[{\rm (i)}]
$\Es_1 W$ has a basis 
$$
  w.
$$
\item[{\rm (ii)}]
For $2 \le i \le D-1$ the subspace $\Es_i W$ has a basis 
$$  
  \Es_i A_{i-1} w, \qquad \qquad \qquad 
  \Es_i A_{i+1} w.
$$
\item[{\rm (iii)}]
$\Es_D W$ has a basis 
$$
  \Es_D A_{D-1} w.
$$
\end{itemize}
\end{corollary}

\begin{corollary}
\label{nime:cor4}
With reference to Notation \ref{blank} let $W$ denote an irreducible $T$-module with endpoint $1$
and local eigenvalue $a_1$. Then the following {\rm (i)--(iii)} hold.
\begin{itemize}
\item[{\rm (i)}]
$\Es_1 W$ has dimension $1$.
\item[{\rm (ii)}]
$\Es_i W$ has dimension $2$ for $2 \le i \le D-1$.
\item[{\rm (iii)}]
$\Es_D W$ has dimension $1$.
\end{itemize}
\end{corollary}
\proof Immediate from Corollary \ref{nime:cor3}. \qed

\begin{corollary}
\label{nime:cor5}
With reference to Notation \ref{blank} let $W$ denote an irreducible $T$-module with endpoint $1$
and local eigenvalue $a_1$. Then the following {\rm (i), (ii)} hold.
\begin{itemize}
\item[{\rm (i)}]
The dimension of $W$ is $2D-2$.
\item[{\rm (ii)}]
The diameter of $W$ is $D-1$.
\end{itemize}
\end{corollary}
\proof Immediate from Corollary \ref{nime:cor4}. \qed

\begin{corollary}
\label{nime:cor6}
With reference to Notation \ref{blank} let $W$ denote an irreducible $T$-module with endpoint $1$
and local eigenvalue $a_1$. Then $W=M^*M w$ for all nonzero $w \in \Es_1 W$.
\end{corollary}
\proof By construction $M^*M w \subseteq W$ and equality holds in view of Theorem \ref{nime:thm2}. \qed


\section{The irreducible $\boldsymbol{T}$-modules with endpoint $\boldsymbol{1}$ 
         and local eigenvalue $\boldsymbol{a_1}$: the $\boldsymbol{A}$-action}
\label{sec:nimeact}

With reference to Notation \ref{blank} let $W$ denote irreducible $T$-module with endpoint 1
and local eigenvalue $a_1$.
In this section we display the action of $A$ on the basis for $W$ given in Theorem \ref{nime:thm2}.
Since $A=L+F+R$ it suffices to give the actions of $L$, $F$, $R$ on this basis.

\begin{lemma}
\label{nimeact:lem1}
With reference to Notation \ref{blank} let $W$ denote an irreducible $T$-module with endpoint $1$
and local eigenvalue $a_1$. Then the following {\rm (i)--(v)} hold for all nonzero $w \in \Es_1 W$.
\begin{itemize}
\item[{\rm (i)}]
$L w = 0$.
\item[{\rm (ii)}]
$L \Es_2 A w = (k - c_2(a_1+1)) w$.
\item[{\rm (iii)}]
For $3 \le i \le D$,
$$
  L \Es_i A_{i-1} w = b_{i-1} \Es_{i-1} A_{i-2} w + (c_i - c_{i-1}) \Es_{i-1} A_i w.
$$
\item[{\rm (iv)}]
$L \Es_2 A_3 w = - b_2 (a_1+1) w$.
\item[{\rm (v)}]
For $3 \le i \le D-1$,
$$
  L \Es_i A_{i+1} w = b_i \Es_{i-1} A_i w.
$$
\end{itemize}
\end{lemma}
\proof
For each equation of Lemma \ref{prod:lem1}, apply each side to $w$ and simplify using $\Es_1 w = w$
and Corollary \ref{sp:cor1a}. \qed

\begin{lemma}
\label{nimeact:lem2}
With reference to Notation \ref{blank} let $W$ denote an irreducible $T$-module with endpoint $1$
and local eigenvalue $a_1$. Then the following {\rm (i)--(iv)} hold for all nonzero $w \in \Es_1 W$.
\begin{itemize}
\item[{\rm (i)}]
$F w = a_1 w$.
\item[{\rm (ii)}]
For $2 \le i \le D-1$,
\begin{eqnarray*}
F \Es_i A_{i-1} w & = & \big( a_{i-1} + a_1(c_i-c_{i-1}) - c_i (a_1+1)(b^i-b^{i-2}) (b^i-1)^{-1} \big) \Es_i A_{i-1} w \\ 
                  &   & - \, c_i (b^i-b^{i-2}) (b^i-1)^{-1} \Es_i A_{i+1} w.
\end{eqnarray*}                  
\item[{\rm (iii)}]
$F \Es_D A_{D-1} w = \big( a_{D-1} + a_1 (c_D-c_{D-1}) - c_D (a_1+1)(b^D-b^{D-2}) (b^D-1)^{-1} \big) \Es_D A_{D-1} w$.
\item[{\rm (iv)}]
For $2 \le i \le D-1$,
$$
  F \Es_i A_{i+1} w = a_i \Es_i A_{i+1} w.
$$
\end{itemize}
\end{lemma}
\proof
For each equation of Lemma \ref{prod:lem2}, apply each side to $w$ and simplify using $\Es_1 w = w$ ,
$\Es_1 A w = a_1 w$ and Lemma \ref{sp:lem1}. \qed

\begin{lemma}
\label{nimeact:lem3}
With reference to Notation \ref{blank} let $W$ denote an irreducible $T$-module with endpoint $1$
and local eigenvalue $a_1$. Then the following {\rm (i)--(iv)} hold for all nonzero $w \in \Es_1 W$.
\begin{itemize}
\item[{\rm (i)}]
For $1 \le i \le D-1$,
$$
  R \Es_i A_{i-1} w = c_i \Es_{i+1} A_i w.
$$
\item[{\rm (ii)}]
$R \Es_D A_{D-1} w = 0$.
\item[{\rm (iii)}]
For $2 \le i \le D-2$,
\begin{eqnarray*}
R \Es_i A_{i+1} w & = & (a_1+1)(c_{i+1} (b^{i-1}-1) (b^{i+1}-1)^{-1} - c_i) \Es_{i+1} A_i w \\
                  &   & + \, c_{i+1} (b^{i-1}-1) (b^{i+1}-1)^{-1} \Es_{i+1} A_{i+2} w.
\end{eqnarray*}                  
\item[{\rm (iv)}]
$R \Es_{D-1} A_D w = (a_1 + 1) (c_D (b^{D-2}-1) (b^D-1)^{-1} - c_{D-1}) \Es_D A_{D-1} w$.
\end{itemize}
\end{lemma}
\proof
For each equation of Lemma \ref{prod:lem3}, apply each side to $w$ and simplify using $\Es_1 w = w$,
$\Es_1 A w = a_1 w$ and Lemma \ref{sp:lem1}. \qed

\section{The isomorphism class of an irreducible $\boldsymbol{T}$-module
         with endpoint $\boldsymbol{1}$}
\label{sec:iso}

With reference to Notation \ref{blank}, in this section we prove
that up to isomorphism there exist exactly two irreducible $T$-modules with endpoint $1$.

\begin{proposition}
\label{iso:prop1}
With reference to Notation \ref{blank}, any two irreducible $T$-modules with endpoint $1$
and local eigenvalue $-1$ are isomorphic.
\end{proposition}
\proof Let $W$ and $W'$ denote irreducible $T$-modules with endpoint $1$ and local eigenvalue $-1$.
Fix nonzero $w \in \Es_1 W$, $w' \in \Es_1 W'$. By Theorem \ref{jeme:thm5},
$W$ and $W'$ have bases $\{\Es_i A_{i-1} w  \mid 1 \le i \le D\}$ 
and $\{\Es_i A_{i-1} w' \mid 1 \le i \le D\}$, respectively.
Let $\sigma \colon W \to W'$ denote the vector space isomorphism defined by 
$\sigma(\Es_i A_{i-1} w) = \Es_i A_{i-1} w'$ for $1 \le i \le D$. 
We show that $\sigma$ is a $T$-module isomorphism. Since $A$ generates $M$ and 
$\Es_0, \Es_1, \ldots, \Es_D$ is a basis for $M^*$, it suffices to show that $\sigma$ commutes with 
each of $A, \Es_0, \Es_1, \ldots, \Es_D$.

\noindent
Using the assertion (iv) below the line \eqref{den0} and the definition of $\sigma$ we immediately 
find that $\sigma$ commutes with each of $\Es_0, \Es_1, \ldots, \Es_D$. It follows from Lemmas 
\ref{jemeact:lem1}--\ref{jemeact:lem3} that $\sigma$ commutes with each of 
$L$, $F$, $R$. Recall $A=L+F+R$ so $\sigma$ commutes with $A$. The result follows. \qed

\begin{proposition}
\label{iso:prop4}
With reference to Notation \ref{blank}, any two irreducible $T$-modules with endpoint $1$
and local eigenvalue $a_1$ are isomorphic.
\end{proposition}
\proof 
Similar to the proof of Proposition \ref{iso:prop1}. \qed

\begin{corollary}
\label{iso:cor3}
With reference to Notation \ref{blank} fix a nonzero $w \in \Es_1 V$ which is orthogonal to 
$\sum_{y \in \G(x)} \h{y}$. Assume that $w$ is an eigenvector for $\Es_1 A \Es_1$.
Then $M^* M w$ is an irreducible $T$-module with endpoint $1$.
\end{corollary}
\proof
Let $H$ denote the subspace of $V$ spanned by the irreducible $T$-modules with endpoint $1$.
By construction and Lemma \ref{sp:lem0} $\Es_1 H$ is the orthogonal complement of 
$\sum_{y \in \G(x)} \h{y}$ in $\Es_1 V$. Hence $w \in \Es_1 H$. Note that $Tw \subseteq H$ so
$Tw$ is the orthogonal direct sum of some irreducible $T$-modules of endpoint $1$. Call these 
$T$-modules $W_1, W_2, \ldots, W_s$. We show $s=1$.
By construction and since $w \in \Es_1 V$ there exist $w_i \in \Es_1 W_i \; (1 \le i \le s)$
such that 
\begin{equation}
\label{eq:sum}
w = w_1 + w_2 + \cdots + w_s.
\end{equation}
For $1 \le i \le s$ we have $w_i \ne 0$; otherwise 
$Tw \subseteq W_1 + \cdots + W_{i-1} + W_{i+1} + \cdots W_s$.
We claim that the $T$-modules $W_1, W_2, \ldots, W_s$ are mutually isomorphic.
To see this, recall that $w$ is an eigenvector for $\Es_1 A \Es_1$; let 
$\eta$ denote the corresponding eigenvalue. Applying $\Es_1 A \Es_1$ to each
term in \eqref{eq:sum} and using $\Es_1 A \Es_1 \in T$ we find 
$\Es_1 A \Es_1 w_i = \eta w_i$ for $1 \le i \le s$. Therefore each of 
$W_1, W_2, \ldots, W_s$ has local eigenvalue $\eta$, so $W_1, W_2, \ldots, W_s$ 
are mutually isomorphic by Propositions \ref{iso:prop1}, \ref{iso:prop4}.
We have proved the claim.
We can now easily show that $s=1$. Suppose $s \ge 2$. By construction there
exists $t \in T$ such that $t w = w_1$. We have 
$w_1 = tw = t w_1 + \cdots + t w_s$ and $t w_i \in W_i$ for $1 \le i \le s$. 
Therefore $t w_i = 0$ for $2 \le i \le s$. Now $(t-I)\Es_1$ is zero on $W_1$ and 
nonzero on $W_i$ for $2 \le i \le s$; this contradicts the fact that 
$W_1, W_2, \ldots, W_s$ are mutually isomorphic. We conclude $s=1$. Now $Tw = W_1$
is an irreducible $T$-module with endpoint $1$. 
The result follows since $Tw = M^*Mw$ by Corollary \ref{jeme:cor8} and Corollary \ref{nime:cor6}. \qed

\medskip \noindent
With reference to Notation \ref{blank} recall that $V$ is an orthogonal direct sum of irreducible 
$T$-modules. Let $W$ denote an irreducible $T$-module. By the {\em multiplicity with which} $W$ 
{\em appears in} $V$ we mean the number of irreducible $T$-modules in this sum which are isomorphic 
to $W$. For example the primary module $M \h{x}$ appears in $V$ with multiplicity $1$. 

\begin{theorem}
\label{iso:thm}
With reference to Notation \ref{blank}, up to isomorphism there exist exactly two irreducible 
$T$-modules with endpoint $1$. The first has local eigenvalue $-1$ and appears in $V$ with multiplicity
$k a_1 (a_1+1)^{-1}$. The second has local eigenvalue $a_1$ and appears in $V$ with multiplicity
$b_1 (a_1+1)^{-1}$. 
\end{theorem}
\proof
By Corollary \ref{local1} each irreducible $T$-module with endpoint $1$ has local eigenvalue 
$-1$ or $a_1$. By Proposition \ref{iso:prop1} (resp. Proposition \ref{iso:prop4}) any 
two irreducible $T$-modules with endpoint $1$ and local eigenvalue $-1$ (resp. $a_1$) 
are isomorphic. For $\eta \in \{a_1, -1\}$ let $\mu_{\eta}$ denote the multiplicity with 
which an irreducible $T$-module with endpoint $1$ and local eigenvalue $\eta$ appears in $V$. 
We show that $\mu_{\eta} = k a_1 (a_1+1)^{-1}$ if $\eta=-1$ and $\mu_{\eta} = b_1 (a_1+1)^{-1}$ 
if $\eta=a_1$. Let $H_{\eta}$ denote the subspace of $V$ spanned by all the irreducible $T$-modules 
with endpoint $1$ and local eigenvalue $\eta$. We claim that $\mu_{\eta}$ is equal to the dimension
of $\Es_1 H_{\eta}$. Observe that $H_{\eta}$ is a $T$-module so it is an orthogonal direct
sum of irreducible $T$-modules:
\begin{equation}
\label{ortsum}
H_{\eta} = W_1 + W_2 + \cdots + W_m \qquad \qquad {\rm (orthogonal\ direct\ sum}),
\end{equation}
where $m$ is a nonnegative integer, and where $W_1, W_2, \ldots, W_m$ are irreducible $T$-modules
with endpoint $1$ and local eigenvalue $\eta$. Apparently $m$ is equal to $\mu_{\eta}$. We show $m$
is equal to the dimension of $\Es_1 H_{\eta}$. Applying $\Es_1$ to \eqref{ortsum} we find
\begin{equation}
\label{ortsum1}
\Es_1 H_{\eta} = \Es_1 W_1 + \Es_1 W_2 + \cdots + \Es_1 W_m \qquad \qquad {\rm (direct\ sum}).
\end{equation}
Note that each summand on the right in \eqref{ortsum1} has dimension 1. It follows that $m$
is equal to the dimension of $\Es_1 H_{\eta}$ and the claim is proven. 
Recall that $\sum_{y \in \G(x)} \h{y}$ is an eigenvector for $\Es_1 A \Es_1$ with eigenvalue $a_1$.
Let $U_{\eta}$ denote the set of those vectors in $\Es_1 V$ that are eigenvectors for $\Es_1 A \Es_1$
with eigenvalue $\eta$ and that are orthogonal to $\sum_{y \in \G(x)} \h{y}$. 
By Corollary \ref{local}(ii) the dimension of $U_{\eta}$ is $k a_1 (a_1+1)^{-1}$ 
if $\eta=-1$ and $ k (a_1+1)^{-1} - 1 = b_1 (a_1+1)^{-1}$ if $\eta=a_1$. 
We now show $\Es_1 H_{\eta} = U_{\eta}$.
By \eqref{ortsum1} and Lemma \ref{sp:lem0} we find $\Es_1 H_{\eta} \subseteq U_{\eta}$. 
Pick a nonzero $w \in U_{\eta}$.
By Corollary \ref{iso:cor3} and definition of $H_{\eta}$ we find 
$w \in \Es_1 M^* M w \subseteq \Es_1 H_{\eta}$ implying $U_{\eta} \subseteq \Es_1 H_{\eta}$. 
It follows $U_{\eta} = \Es_1 H_{\eta}$. The result now follows from these comments and the above claim. \qed

\bigskip \noindent
{\bf \large Acknowledgement:} The author would like to thank Paul Terwilliger for his careful 
reading of the earlier versions of this manuscript and many helpful suggestions.

\end{document}